\documentclass[12pt,english,twoside,a4paper]{article}
\usepackage{latexsym}
\usepackage{babel}
\usepackage{amsmath}
\usepackage{amsfonts}
\usepackage{amssymb}
\usepackage{graphicx}
\usepackage{xcolor}
\usepackage{enumitem}
\usepackage{pifont}  
\usepackage{euscript}  
\DeclareFontFamily{OT1}{pzc}{}
\DeclareFontShape{OT1}{pzc}{m}{it}{<-> s * [1.200] pzcmi7t}{}
\DeclareMathAlphabet{\mathpzc}{OT1}{pzc}{m}{it}
\newcommand{\enne}{\mathbb{N}}
\newcommand{\erre}{\mathbb{R}}

\newcommand{\bg}{\bigskip}
\newcommand{\md}{\medskip}
\newcommand{\sm}{\smallskip}
\newcommand{\np}{\noindent}
\newcommand{\beq}{\begin{equation}}
\newcommand{\eeq}{\end{equation}}
\newcommand{\oP}{\overline{P_\rho}}

\newcommand{\les}{\preccurlyeq}
\newcommand{\ges}{\succcurlyeq}

\newcommand{\E}{\mathpzc{E}}
\newtheorem{theorem}{Theorem}[section]
\newtheorem{corollary}{Corollary}[section]

\newtheorem{definition}{Definition}[section]
\newtheorem{remark}{Remark}[section]

\newtheorem{example}{Example}[section]
\title{A positive eigenvalue result for semilinear differential equations in Banach spaces with functional initial conditions}
\author{{\bf  Gennaro Infante and Paola Rubbioni}\\
{\small {\em Department of Mathematics and Computer Science, University of Calabria, Italy}}\\ {\small {\em Department of Mathematics and Computer Science, University of Perugia, Italy}}\\
{\small {\em E-mail addresses: gennaro.infante@unical.it,  paola.rubbioni@unipg.it}}
}
\date{}
\begin{document}

\maketitle

\np 
{\bf Abstract.} We study the existence of positive eigenvalues with associated nonnegative mild eigenfunctions for a class of abstract initial value problems in Banach spaces with functional, possibly nonlocal, initial conditions. The framework includes periodic, multipoint, and integral average conditions. Our approach relies on nonlinear analysis, topological methods, and the theory of strongly continuous semigroups, yielding results applicable to a wide range of models. As an illustration, we apply the abstract theory to a reaction–diffusion equation with a nonlocal initial condition arising from a heat flow problem.
\md

\np
\md

\np

\bg
\bg

\section{Introduction}
In this manuscript we study the existence of \emph{positive} eigenvalues $\lambda$, with associated nonnegative eigenfunctions $y$, of the initial value problem 
\begin{equation}\label{e:P}
\begin{cases}
y'(t)= Ay(t)+\lambda f(t,y(t)) ,\ t \in [0,1],\\
y(0)=\lambda B(y),
\end{cases}
\end{equation}
where $A$ is a linear operator defined on a dense subset of a Banach space $E$ and taking values in $E$ itself, $f:[0,1]\times E \to E$ and $B:C([0,1];E)\to E$ are given functions.

We stress that the eigenfunctions are understood in the mild sense and that the initial condition is of functional type, not necessarily linear or local.
There exists a wide body of literature on the solvability of the initial value problem \eqref{e:P} in the case $\lambda = 1$, under a variety of initial conditions, both for differential equations and for differential inclusions. We refer, for instance, to the works \cite{b:b-t-v, b:Bou1, b:Bou3, CPR15, MR16, O13} and to the books \cite{KOZ, V03}.
The situation we consider here—namely, an initial datum of functional type—is fairly broad, encompassing not only periodic conditions but also nonlocal initial conditions, such as the multipoint condition introduced by Byszewski in \cite{B91} or the integral average condition used by several authors (see, for example, \cite{MR16} or \cite{YuWa22}). In particular, our interest in the presence of functional conditions in \eqref{e:P} stems from the fact that they can be used to model physical phenomena. To illustrate this, just think about the  parabolic boundary value problem
\beq
\label{e-Ex-intro}
\begin{cases}
u_t(t,x)= u_{xx}(t,x)+ \lambda\, g(t,x,u(t,x)), \quad (t,x)\in (0,1)\times (0,\pi),\\
u(0,x)= \lambda \, \alpha(x)\int_0^1 u(s,\pi/2)\, ds, \quad x\in (0,\pi),\\
u(t,0)=u(t,\pi)=0, \quad t\in (0,1).
\end{cases}
\eeq
This is a model involving a one-dimensional heat equation subject to an external source, where a sensor is placed at the midpoint of a thin rod of length $\pi$ and one seeks solutions that, when evaluated at time $0$, satisfy a time-averaged condition on the temperature measured at the midpoint of the bar.

Our approach is fairly general and relies on tools from nonlinear analysis in abstract spaces, topological methods, and the theory of strongly continuous semigroups; as a consequence, the results we obtain are applicable not only to problem \eqref{e-Ex-intro}, but also to all models that can be rewritten in the form \eqref{e:P}.
In order to make the paper as self-contained as possible, Section \ref{s:Prel} collects the necessary definitions and preliminary results that will be repeatedly used in the subsequent sections. Then, in Section \ref{s:Prob} we present the framework of the problem and state the main result of the paper. In the subsequent section, namely Section \ref{s-Appl}, we show how our general result can be applied to the case of a reaction–diffusion equation, and we also provide an example of the form \eqref{e-Ex-intro}.

\section{Preliminaries}
\label{s:Prel}

We recall the next definitions (see, e.g., \cite[Sections 8.1 and 8.2]{AB06}).
\sm

\np
A set $C$ of a vector space $X$ is a {\em pointed convex cone} if 
\begin{description}
\item{(C1)} $C+C \subset C$;
\item{(C2)}  $\alpha\, C\subset C$ for all $\alpha\ge 0$;
\item{(C3)}  $C\cap (-C)=\{0\}$. 
\end{description}
A pointed convex cone $C$ induces on $X$ a partial ordering $\les$  by
\[
x\les y \ \Leftrightarrow \ y-x\in C.
\]
This partial ordering is compatible with the algebraic structure of the space, that is the next properties hold: 
\begin{description}
\item{(1)} $x\les y \ \Rightarrow \ x+z \les y+z, \mbox{ for every } x,y,z\in X;$
\item{(2)} $x\les y \ \Rightarrow \ \alpha x \les \alpha y, \mbox{ for every } x,y\in X,\, \alpha \ge 0.$
\end{description}

\np
An {\em ordered vector space} $(X, \les)$  is a real vector space with a partial order relation $\les$ that is compatible with the algebraic structure of the space, that is (1) and (2) are satisfied. In this case, the set 
\[
X^+:=\{x\in X: 0\les x\}
\]
is a pointed convex cone, called the {positive cone} of $X$. On the other hand, any vector space $X$ endowed with a pointed convex cone is an ordered vector space.

\np
A {\em vector lattice} $(X, \les)$ is an ordered vector space such that every pair of elements has a {\em supremum} and an {\em infimum}, where an element $z$ is the supremum of a pair $x,y\in X$ and we write $z:=x\vee y$ if
\begin{description}
\item{(i)} $x\les z$ and $y\les z$;
\item{(ii)} if $x\les u$ and $y\les u$, then $z\les u$.
\end{description}
The definition of infimum is analogous.
Moreover, in a vector lattice the {\em positive part} and {\em negative part} of any element $x\in X $ exist, i.e. \[
x^+=x\vee 0\ \mbox{ and } \  x^-=-x\vee 0.
\]
Moreover, it is defined the {\em absolute value} of $x$ as
\[
|x|:=x\vee -x.
\]
Of course, it holds that
\[x=x^+ - x^-  \ \mbox{ and } \ |x|=x^+ + x^- .\]

\np A {\em normed vector lattice} $(X, \les, \|\cdot\|)$ is a  vector lattice where the norm is such that
\beq\label{e-an}
|x|\les |y| \ \Rightarrow \ \|x\|\le \|y\|.
\eeq

A  {\em Banach lattice} $(X, \les, \|\cdot\|)$ is a complete normed vector lattice  (see, e.g. \cite[Section 9.1]{AB06}). The Euclidean spaces $\erre^n$ with their Euclidean norms are all Banach lattices; also, if $K$ is a compact space, then the space $C(K)$ of all the real continuous functions over $K$ under the sup norm is a Banach lattice. For further examples of Banach lattices we refer to \cite[Example 9.1]{AB06}.
\md

Let $(\Xi, \Sigma, \mu)$ be a measure space and \((X,\|\cdot\|_X) \) be a real Banach space. 
\\
A function $f:\Xi\to X$ is {\em strongly $\mu$-measurable} if there exists a sequence $(\varphi_n)_n$ of $X$-simple functions such that $\lim_{n\to \infty} \|f(\xi)-\varphi_n(\xi)\|_X=0$ for $\mu$-almost all $\xi\in\Xi$. 
\\
A strongly  $\mu$-measurable function $f:\Xi\to X$ is {\em Bochner integrable} if there exists a sequence $(\varphi_n)_n$ of $X$-step functions such that the real measurable function $\|f-\varphi_n\|_X$ is Lebesgue integrable for every $n\in\enne$ and $\lim_{n\to\infty}\int \|f-\varphi_n\|_X\, d\mu=0$. Then, for each $A\in \Sigma$ the {\em Bochner integral of $f$ over $A$} is 
\[\int_A f\, d\mu:=\lim_{n\to\infty} \int_A \varphi_n\, d\mu,\]
where the limit is in the norm topology on $X$. 
\\
Moreover, if $X$ is a Banach lattice and $f,g:\Xi \to X$ are Bochner integrable functions such that $f(\xi)\les g(\xi)$ for $\mu$-almost all $\xi\in\Xi$, then
\[\int_A f\, d\mu\les \int_A g\, d\mu,\]
for every $A\in\Sigma$ (cf. \cite[Theorem 11.43]{AB06}). For further properties of the Bochner integral we refer, e.g., to \cite{AB06}.
\sm

From now on, we denote by
$L^1([a,b];X)$ the set of the functions which are Bochner integrable on $[a,b]$ with respect to the Lebesge measure, and by $C([a,b];X)$ its subset of continuous functions ($C([a,b])$ if $X=\erre$). In $C([a,b];X)$ with the usual sup-norm that we denote by $\|\cdot\|_{C([a,b];X)}$. Moreover, by $L^p([a,b];X)$ we denote the set of functions such that their $p$-power is Bochner integrable on $[a,b].$
\md

Let ${\cal K}$ be a {\em cone} in a Banach space $(X,\|\cdot\|)$, that is a closed pointed convex cone. 
In our main result we will make use of the following 
Birkhoff–Kellogg type theorem due Krasnosel'ski\u{\i} and Lady\v{z}enski\u{\i}~\cite{Kra-Lady}.
\begin{theorem}\cite[Theorem~5.5]{Krasno}
\label{t-BK}
Let \( (X, \| \cdot \|) \) be a real Banach space, \( U \subset X \) an open bounded set with \( 0 \in U \), \( {\cal K} \subset X \) a cone, and let \( T \colon {\cal K} \cap \bar U \to {\cal K} \) be a compact operator (i.e. it is continuous and maps bounded sets into relatively compact sets). Suppose that
\[
\inf_{x \in {\cal K} \cap \partial U} \|Tx\| > 0.
\]
Then there exist \( \lambda_0 \in (0, +\infty) \) and \( x_0 \in {\cal K} \cap \partial U \) such that
\( 
x_0 = \lambda_0 T x_0. 
\)
\end{theorem}

\np
Moreover, let us recall that a cone ${\cal K}$ is said to be  {\em normal} \  if 
there exists $c >0$ such that
\beq
\label{e-Knorm}
\|x\| \le c\, \|y\|,\ \mbox{ for every $x,y\in K$ with $x\les y$, }
\eeq
where $\les$ is the partial ordering induced by ${\cal K}$.

\begin{remark}\label{r-1}
We observe that if $(X,\les, \|\cdot\|)$ is a Banach lattice, then the positive cone $X^+$ of $X$ is a normal cone. Indeed, let $x,y\in X^+$ with $x\les y$ be fixed. Of course, $0\les x$ and $ 0\les y$, so that $x=|x|$, $y=|y|$. Thus, by $x\les y$ we have $|x|\les |y|$. By using now the compatibility of the norm with the ordering (see \eqref{e-an}), we get $\|x\|\le \|y\|$, that is \eqref{e-Knorm} for $c=1$.
\end{remark}

Let ${\cal L}(X)$ be  the Banach space of all bounded linear operators from the Banach space $(X,\|\cdot\|)$  to itself, furnished with the usual operator norm $\|L\|_{{\cal L}}=\sup_{\|x\|\le 1}\|L x\|$. 
A family $\{U(t)\}_{t \geq 0} $ of bounded linear operators on the Banach space $X$  is called a \emph{$C_0$-semigroup} (or a \emph{strongly continuous semigroup}) if it satisfies the following properties (see, e.g. \cite[Definition I.1.1]{EN06}):
\begin{description}
    \item{(U1)}  $U(0) = I$, the identity operator on $X$;
     \item{(U2)}   $U(t+s) = U(t)U(s)$ for all $t, s \geq 0$;
      \item{(U3)}  for every $x\in X$, the orbit maps $\xi_x:t\mapsto \xi_x(t):=U(t)x$ are continuous from $\erre^+$ to $X$.
\end{description}
It is known that (cf. \cite[Proposition I.1.4]{EN06}) 
for every $C_0$-semigroup $\{U(t)\}_{t \geq 0}$  there exist $\hat{M}\ge 1$ and $\delta\in \erre$ such that
\begin{equation}\label{e:Mdelta}
\|U(t)\|_{{\cal L}}\le \hat{M} e^{\delta t}, \mbox{  for every $t\ge 0$}.
\end{equation}
Further (see, e.g. \cite[Definition 2.3.1]{P83}) a $C_0$-semigroup \( \{U(t)\}_{t \geq 0} \)  is said to be {\em compact}   if 
 $U(t)$  is a compact operator, for every $t > 0$ (i.e., each $U(t)$ maps bounded sets into relatively compact sets).

\label{t-Uc}\begin{theorem}\cite[Theorem 2.3.2]{P83} Let $\{U(t)\}_{t\ge 0}$ be a $C_0$ semigroup. If $U(t)$ is compact for $t > t_0$, then $U(\cdot)$ is continuous in the uniform operator topology for $t > t_0$.
\end{theorem}

\np
If  $(X, \ges, \|\cdot\|)$ is a Banach lattice, we say that a \( C_0 \)-semigroup \( \{U(t)\}_{t \ge 0} \) is {\em positive} if 
\beq
\label{e-Up} 
U(t) x  \ges 0  \mbox{ for all } t \ge 0 , \mbox{ whenever }  x\ges 0.
\eeq 
\sm

\np
An operator  \( A : D(A)\subset X \to X \)  is said to be the {\em (infinitesimal) generator of a \( C_0 \)-semigroup \( \{U(t)\}_{t\ge 0} \) } if (see, e.g. \cite[Definition II.1.2]{EN06})
\[
Ax := \dot\xi_x(0)=\lim_{t \downarrow 0} \frac{1}{t} \left( U(t) - \operatorname{id} \right)x,
\]
for every $x$ in its domain
\[D(A) :=\left\{ x \in X : \xi_x \mbox{  is differentiable in } \erre^+ \right\}.\]
Now,  for a strongly continuous semigroup \( \{U(t)\}_{t \geq 0} \) and an element $x\in X$, each orbit map $\xi_x:t\mapsto U(t)x$ is differentiable in $\erre^+$ if and only if $\xi_x$ is right differentiable at $t=0$ (cf. \cite[Lemma II.1.1]{EN06}). Hence, the domain of $A$ can be written as
\[
D(A) = \left\{ x \in X : \lim_{t \downarrow 0} \frac{1}{t} \left( U(t) - \operatorname{id} \right)x \text{ exists in } X \right\}.
\]
For the generator $A$ of a $C_0$-semigroup  \( \{U(t)\}_{t \geq 0} \) the following properties hold:
\begin{description}
\item{(A1)} $A$ is a linear operator;
\item{(A2)} if $x\in D(A)$, then for every $t\ge 0$ it holds that  $U(t)x\in D(A)$ and 
\[\frac{d}{dt}U(t)x=U(t)Ax=AU(t)x;\]
\item{(A3)} for every $t\ge 0$ and $x\in X$ one has
\[\int_0^t U(s)x\, ds \in D(A).\]
\end{description}

Finally, in order to prove our main result, we will need the Ascoli-Arzelà theorem in Banach spaces.

\begin{theorem}\cite[Theorem A.2.1]{V03} 
\label{t-AA}
A subset $\mathcal{F}$ in $C([a,b]; X)$ is relatively compact if and only if:
\begin{enumerate}
    \item $\mathcal{F}$ is equicontinuous on $[a,b]$;
    \item there exists a dense subset ${\cal D}$ in $[a,b]$ such that, for each $t \in {\cal D}$, $\mathcal{F}(t) = \{f(t); f \in \mathcal{F}\}$ is relatively compact in $X$.
\end{enumerate}
\end{theorem}

\section{Problem setting and main result}
\label{s:Prob}

In this section, we provide (in the mild sense) eigenvalues and eigenfunctions of the initial value problem with functional initial conditions \eqref{e:P}.
\sm

Let $(E,\ges,\|\cdot\|)$ be a Banach lattice, where the partial ordering $\ges$ is induced by a given normal cone $K$. We denote by $P$ the corresponding positive cone in $C([0,1];E)$, i.e. 
\[ 
P=\{y\in C([0,1];E):  y(t)\ges 0, \mbox{ for every } t\in [0,1]\}. 
\]
Since $E$ is a Banach lattice, the space $C([0,1];E)$ with the partial ordering induced by $P$ is a Banach lattice too. 
For every $\rho\in (0,+\infty)$, we consider the sets
%
\[ K_\rho:= \{v\in K: \|v\|< \rho\}, \quad \overline{K_\rho}:= \{v\in K: \|v\|\le  \rho\}, \quad \partial K_\rho:= \{v\in K: \|v\|= \rho\}, \]
and 
\[ P_\rho:= \{y\in P: \|y\|_C< \rho\}, \quad \overline{P_\rho}:= \{y\in P: \|y\|_C\le \rho\}, \quad \partial P_\rho:= \{y\in P: \|y\|_C= \rho\}. \]

Suppose that 
$A$ is the infinitesimal generator of a $C_0$-semigroup $\{U(t)\}_{t\ge 0}$.
\sm

\begin{definition}
We say that a solution of \eqref{e:P} is a couple $(\lambda, y)$, where $\lambda\in \mathbb{R}$ and $y:[0,1]\to E$ is a continuous function, such that the integral equation
\begin{equation*}
    y(t)=\lambda\, U(t)B(y)+\lambda\, \int_0^t U(t-s) f(s,y(s))\, ds, \ t\in [0,1],
\end{equation*}
is satisfied.
Then we say that $y$
is a {\it (mild) eigenfunction} of~\eqref{e:P} with a corresponding eigenvalue $\lambda$.
\end{definition}

\begin{theorem}\label{T-main}
Let $\rho\in (0,+\infty)$. Assume the following conditions hold.
\begin{description}
\item{(H1)} $A$ is the infinitesimal generator of a positive and compact $C_0$-semigroup $\{U(t)\}_{t\ge 0}$;
\item{(H2)} $f:[0,1]\times \overline{K}_\rho \to K$ 
is a Carathèodory function 
mapping bounded sets into bounded sets
and there exists $\delta_\rho:[0,1]\to K$ such that
\beq
\label{e-fdelta} 
f(t,y(t)) \ges \delta_\rho(t), \ \mbox{for every } t\in [0,1], y\in  \partial P_\rho;
\eeq
\item{(H3)} $B: C([0,1];\overline{K}_\rho) \to K$ 
 is a continuous function mapping bounded sets into bounded sets 
and there exists $\eta_\rho\in K$ such that 
\beq
\label{e-Beta}
 B(y)\ges \eta_\rho, \ \mbox{for every } y\in \partial P_\rho;
 \eeq
\item{(H4)} there exists $t_0\in\, ]0,1]$
such that
\beq
\label{e-H4} 
\Bigl\|U(t_0)\eta_\rho + \int_0^{t_0} U(t_0-s) \delta_\rho(s)\, ds\Bigr\|>0,
\eeq
where \( \{U(t)\}_{t \geq 0} \) is the $C_0$-semigroup generated by $A$ and $\eta_\rho, \, \delta_\rho$ are from (H3), (H2), respectively.
\end{description}
Then there exists a positive solution of \eqref{e:P}, i.e. an eigenvalue $\lambda_\rho\in (0,+\infty)$ and a (mild) eigenfunction $y_\rho\in P$, such that $\|y_\rho\|_C=\rho$.
\end{theorem}

\np
{\bf Proof.} First, recall that every $C_0$-semigroup satisfies the condition \eqref{e:Mdelta} for some $\tilde M\ge 1$ and $\delta \in \erre$. So this also holds for the $C_0$-semigroup  generated by $A$ according to the hypothesis (H1). Hence, we have   
\beq\label{e-D} 
\|U(t)\|_{{\cal L}}\le \bigl(\hat{M} \sup_{t\in [0,1]}e^{\delta t} \bigr)=:D, \mbox{ for every } t\in [0,1].
\eeq

Let us consider the operator $T:\overline{P_\rho} \to P$ defined by
\beq
\label{e-Tdef}
  Ty(t)=U(t)B(y)+\int_0^t U(t-s) f(s,y(s))\, ds, \ t\in [0,1]. 
  \eeq

We first show that the operator $T$ is well-posed. Let $y\in \overline{P_\rho}$ be arbitrarily fixed. First of all, note that the function $Ty$ is continuous on $[0,1]$ by the properties of the $C_0$-semigroups and of the Bochner integral. Moreover, by hypotheses (H2) and (H3) we have 
\[
B(y)\ges 0 \quad \mbox{ and } \quad f(t,y(t)) \ges 0.
\]
Since the $C_0$-semigroup \( \{U(t)\}_{t \geq 0} \) is positive (see (H1)) and the Bochner integral is monotone
with respect to the ordering given by the cone $K$, it holds that $Ty(t)\ges 0$ for every $t\in [0,1]$. Hence, the function $Ty$ belongs to $P$.
\md

Now, we prove that $T$ is a compact operator. 
\\
We begin by showing the continuity of $T$.
Let $y_*$ an arbitrary element of $\overline{P_\rho}$ and $\{y_n\}_n $ be any sequence in $\overline{P_\rho}$ converging to $y_*$ .
For  every $t\in [0,1]$ and $n\in\enne$, we have  
\begin{align*}
\|Ty_n(t)-Ty_*(t)\| 
& = \Bigl\|U(t)B(y_n)+\int_0^t U(t-s) f(s,y_n(s))\, ds +\\
& \qquad\qquad \qquad - U(t)B(y_*)-\int_0^t U(t-s) f(s,y_*(s))\, ds\Bigr\|\\
& \le \|U(t)\|_{{\cal L}}\|B(y_n)- B(y_*)\|+\\
& \qquad\qquad \qquad +\int_0^t \|U(t-s)\|_{{\cal L}}\|f(s,y_n(s))-f(s,y_*(s))\|\, ds \\
 & \le D\|B(y_n)- B(y_*)\|+D\int_0^1 \|f(s,y_n(s))-f(s,y_*(s))\|\, ds,
\end{align*}
where $D$ is from \eqref{e-D}. Hence, for every $n\in \enne $ the following estimate holds
\beq
\label{e:stT}
\|Ty_n-Ty_*\|_C\le D\|B(y_n)- B(y_*)\|+D\int_0^1 \|f(s,y_n(s))-f(s,y_*(s))\|\, ds.
\eeq
The set of functions $\oP$ is bounded, which implies the uniform boundedness of the set 
\[\oP([0,1]):= \bigcup\limits_{s\in [0,1]}\oP(s),\]
where $\oP(s):=\{y(s)\in E : y\in \oP\}$. Note that the set $[0,1]\times \oP([0,1])$ is also bounded. Therefore, condition (H2) implies that $f\left( [0,1]\times \oP([0,1]) \right)$ is a bounded subset of $E$, 
so that 
there exists $M_\rho>0$ for which
\beq
\label{e-Mrho} \|f(s,y(s))\|\le M_\rho, \mbox{ for every } s\in [0,1],\ y\in \oP.
\eeq
Hence we get 
\[ \|f(s,y_n(s))-f(s,y_*(s))\|       %
  \le 2M_\rho,  \mbox{ for every $s\in [0,1],\ n\in\enne.$}
\]
As a consequence, we can apply the Lebesgue dominated convergence Theorem, and by \eqref{e:stT} we obtain 
\[
\lim_{n\to +\infty} \|Ty_n-Ty_*\|_C\le 
D\lim_{n\to +\infty}\|B(y_n)- B(y_*)\|+D\int_0^1 \lim_{n\to +\infty} \|f(s,y_n(s))-f(s,y_*(s))\|\, ds.
\]
By the continuity properties given by (H2) and (H3) we obtain  
\[
\lim_{n\to +\infty} \|Ty_n-Ty_*\|_C=0.
\]
The continuity of $T$ follows from the arbitrariness of $y_*$ and $\{y_n\}_n$.
\sm

\np
We prove now that $T$ maps bounded sets into relatively compact sets. 
\\
To this aim, we first show that $T(\oP)$ is relatively compact. To this aim, we use  Theorem \ref{t-AA}.
Let us split the operator $T$ as 
\beq\label{e-THG}T=H+G,\eeq
where $H:\overline{P_\rho}\to C([0,1];E)$ is the nonlocal operator 
\[Hy(t)=U(t)B(y),\ t\in [0,1], 
\]
and $G:\overline{P_\rho}\to C([0,1];E)$ is the Cauchy operator 
\[
Gy(t)=\int_0^t U(t-s)f(s,y(s))\, ds,\ t\in [0,1].
\]
For every $t\in [0,1]$, each set 
\[H(\oP)(t)\equiv U(t)B(\oP)\]  
is relatively compact. Indeed, the set $\oP$ is bounded, so that, by assumption (H3), the set $B(\oP)$ is also bounded. Then, for every $t\in [0,1]$, we use the compactness of $U(t)$ guaranteed by the hypothesis (H1) and achieve the relative compactness of $H(\oP)(t)$.
\\
We show now the relative compactness of 
\[ G(\oP)(t):= \left\{ \int_0^t U(t-s)f(s,y(s))\, ds : y\in \oP \right\},\]
for every $t\in [0,1]$. 
\\
If $t=0$, note that the set $G(\oP)(0)=\{0\}$ is compact. 
\\
Let us fix $t\in \, ]0,1]$. We consider a real number $\eta$ with $0< \eta < t$ and put
\beq
\label{e-Geta} 
G_\eta y(t):= \int_0^{t-\eta} U(t-s)f(s,y(s))\, ds=U(\eta)\int_0^t U(t-s-\eta)f(s,y(s))\, ds. 
\eeq
The set
\[ G_\eta(\oP)(t):=\{G_\eta y(t): y\in \oP\} \]
is relatively compact. 
In fact, we know by \eqref{e-Mrho} that the set \( \{f(s,y(s)) : s\in [0,1], \, y\in \oP\} \)
is bounded. By the compactness of $U(t)$ for every $t>0$, the set 
\[\left\{U(\eta)\int_0^{t-\eta} U(t-s-\eta)f(s,y(s))\, ds : y\in \oP\right\}\]is relatively compact.
That is the relative compactness of $ \{G_\eta y(t): y\in \oP\} $ (cf. \eqref{e-Geta}).
\\
Therefore, for $y\in \oP$, by using also \eqref{e-D} and \eqref{e-Mrho}, we have the estimate
\begin{align*}
\|Gy(t)-G_\eta y(t)\|
&=\left\|\int_{t-\eta}^t U(t-s) f(s,y(s))\, ds\right\|\le \eta DM_\rho,
\end{align*}
which implies the totally boundedness of $G(\oP)(t)$, and then its relative compactness. 
\\
Finally, by \eqref{e-THG} we achieve the relative compactness of $T(\oP)(t)$, for every $t\in [0,1]$.
\md

\np
Let us prove that $T(\oP)$ is a family of equicontinuous functions. Let $t_1,t_2\in \, ]0,1]$ -- without restriction of generality, assume $t_1<t_2$ -- and $y\in \oP$. By using again \eqref{e-D} and \eqref{e-Mrho}, it holds that
\begin{align*}
\|Ty(t_2)-Ty(t_1)\|
& \leq \| [U(t_2)-U(t_1)]B(y) \| + \left\| \int_{t_1}^{t_2} U(t_2-s)f(s,y(s))ds \right\|+\\
&\quad  
+\left\| \int_{0}^{t_1} [U(t_2-s)-U(t_1-s)]f(s,y(s))ds \right\| \\
&\le \|U(t_2)-U(t_1)\|_{{\cal L}}\|B(y)\| + \int_{t_1}^{t_2} \left\|U(t_2-s) \right\|_{{\cal L}} \left\|f(s,y(s)) \right\|\, ds \, +\\
&\quad + \int_0^{t_1}\left\| U(t_2-s)-U(t_1-s) \right\|_{{\cal L}}  \left\|f(s,y(s)) \right\|\, ds\,  \\
&\le  \|U(t_2)-U(t_1)\|_{{\cal L}}\|B(y)\| +  D M_\rho (t_2-t_1)+ \\
&\quad + M_\rho\int_0^{t_1}\left\| U(t_2-s)-U(t_1-s) \right\|_{{\cal L}}  \, ds\, .
\end{align*} 
Moreover, by (H3) there exists $N_\rho>0$ for which 
\beq
\label{e-Nrho} \|B(y)\|\le N_\rho, \mbox{ for every }  y\in \oP,
\eeq
so that 
\begin{align}
\nonumber
\|Ty(t_2)-Ty(t_1)\|
&\le  N_\rho \|U(t_2)-U(t_1)\|_{{\cal L}} +  D M_\rho (t_2-t_1) +\\
\label{e-Uconvi}
&\quad + M_\rho\int_0^{t_1}\left\| U(t_2-s)-U(t_1-s) \right\|_{{\cal L}}  \, ds\, .
\end{align} 
Now, each operator $U(t)$, $t>0$, is compact (see (H1)), hence Theorem \ref{t-Uc} yields that the operator $t\mapsto U(t)$ is continuous in the uniform operator topology for  $t>0$. Therefore for $t_2\rightarrow t_1$ we get
\[
\|U(t_2)-U(t_1)\|_{{\cal L}} \rightarrow 0 , \quad \mbox{ and } \quad
\|U(t_2-s)-U(t_1-s)\|_{{\cal L}} \rightarrow 0 \mbox{ for every $0\le s \le t_1$.}
\]
Since for every $0\le s\le t_1\le t_2$ it is true that
\[
\|U(t_2-s)-U(t_1-s)\|_{{\cal L}}\le 2D , 
\]
then we can apply the Lebesgue dominated convergence Theorem and, by \eqref{e-Uconvi}, we obtain 
\begin{align*}
0\le \|Ty(t_2)-Ty(t_1)\|
&\le  N_\rho \|U(t_2)-U(t_1)\|_{{\cal L}} +  D M_\rho (t_2-t_1) +\\
&\quad + M_\rho\int_0^{t_1}\left\| U(t_2-s)-U(t_1-s) \right\|_{{\cal L}}  \, ds\xrightarrow[t_2\to t_1]{} 0,
\end{align*} 
independently of $y$.
\\
The functions of the family $T(\oP)$  are also equibounded. Indeed, by \eqref{e-D},  \eqref{e-Nrho}, \eqref{e-Mrho}  for any $y\in\oP$ we have
\[
\|T(y)\|_C
\le \sup_{t\in [0,1]} \left\{ \|U(t)\|_{{\cal L}}\|B(y)\|+\int_0^t \|U(t-s)\|_{{\cal L}} \| f(s,y(s))\|\, ds\right\}
\le D(N_\rho+ M_\rho).
\]
Hence, we can apply the Arzelà-Ascoli Theorem and obtain that $T(\oP)$ is relatively compact.
\sm

\np
Now, let ${\cal D}$ be any bounded subset of $\oP$. Observe that 
\[ \overline{T({\cal D})}  \subset \overline{T(\oP)}, \]
where the set $\overline{T(\oP)}$ is compact, so that $\overline{T({\cal D})}$ is also compact. Thus we obtain the relative compactness of $T({\cal D})$. 
\sm

\np
Therefore, recalling that $T$ is also continuous, we conclude that it is a compact operator.

\bg
In order to apply Theorem \ref{t-BK}, it remains to show that
\beq
\label{e-Tinf}
\inf_{y \in \partial P_\rho } \|T(y)\|_C > 0.
\eeq
To this aim, let us consider any $y\in \partial P_\rho$. 
\\
First of all,  for every fixed $t\in [0,1]$, the following inequality holds:
\beq
\label{e-Umon1}
U(t)B(y)\ges U(t)\eta_\rho,
\eeq
where $\eta_\rho$ is from assumption (H3). In fact, by \eqref{e-Beta} we have $B(y)\ges \eta_\rho $ or, equivalently,
\(
 B(y)- \eta_\rho \ges 0.
 \)
The operator $U(t)$ is positive from (H1), hence   \(U(t)[ B(y)- \eta_\rho] \ges 0\) (cf.~\eqref{e-Up}). Therefore, the linearity of $U(t)$ allows to  the relation
\[
U(t)B(y)- U(t)\eta_\rho \ges 0,
 \]
which is equivalent to \eqref{e-Umon1}.
Analogously, from \eqref{e-fdelta} we get
\beq
\label{e-Umon2}
U(t-s) f(s,y(s)) \ges U(t-s) \delta_\rho(s),\ \mbox{ for every $0\le s\le t$. }
\eeq
Therefore, the conditions \eqref{e-Umon1}, \eqref{e-Umon2}, and the order preserving property
 of the Bochner integral with respect to the cone $K$, imply that
\[
U(t)B(y) + \int_0^t U(t-s)  f(s,y(s)) \, ds
\ges U(t)\eta_\rho + \int_0^t U(t-s) \delta_\rho(s)\, ds.
\]
Since the cone $K$ is normal, we have 
\[
c\left\|U(t)B(y) + \int_0^t U(t-s)  f(s,y(s)) \, ds \right\|
\ge \left\|U(t)\eta_\rho + \int_0^t U(t-s) \delta_\rho(s)\, ds \right\|,
\]
for some $c>0$ (cf. \eqref{e-Knorm}).
\\
Now, let us consider the element $t_0\in \,] 0,1]$  from hypothesis (H4). By the above inequality taken for $t=t_0$ and by using \eqref{e-H4}, we can derive that
\begin{align*}
\|T(y)\|_C
&\ge \|T(y)(t_0)\|=\left\|U(t_0)B(y) + \int_0^{t_0} U(t_0-s)  f(s,y(s)) \, ds \right\|\\
&\ge \dfrac1c \left\|U(t_0)\eta_\rho + \int_0^{t_0} U(t_0-s) \delta_\rho(s)\, ds \right\|.
\end{align*}
Note that the latter is strictly positive by \eqref{e-H4} and does not depend on $y\in \partial P_\rho$. So, 
$
\inf_{y \in \partial P_\rho } \|T(y)\|_C > 0
$, that is \eqref{e-Tinf}.
\sm

Finally, from Theorem \ref{t-BK} it follows that there exist a positive eigenvalue $\lambda_\rho$ and an eigenfunction $y_\rho\in \partial P_\rho$, i.e. 
\[
y_\rho=\lambda_\rho T(y_\rho).
\] 
In other words (cf. \eqref{e-Tdef}), it holds that
 \[
y_\rho=\lambda_\rho U(t)B(y_\rho) +\int_0^t U(t-s)  \lambda_\rho f(s,y_\rho(s)) \, ds.
\] 
The thesis is therefore achieved.
\hfill $\Box$

\begin{corollary}\label{CorEig}
In addition to the hypotheses of Theorem \ref{T-main}, assume that $\rho$  can be chosen arbitrarily in $(0,+\infty)$. Then, for every $\rho$ there exists a non-negative eigenfunction $y_{\rho}\in \partial P_{\rho}$ of the boundary value problem \eqref{e:P} to which corresponds a $\lambda_{\rho} \in (0,+\infty)$.
\end{corollary}

\section{An application to reaction-diffusion equations}\label{s-Appl}
Given $x_*\in \erre^n$ and $R>0$, let $B_R(x_*)\subset\erre^n$ be the open ball with center $x_*$ and radius $R$. Let $\left(C_0(B_R(x_*)), \|\cdot\|_\infty\right)$ be the Banach space  of all real continuous functions defined on $\overline{B}_R(x_*)$ and vanishing on the boundary of $B_R(x_*)$, equipped with the usual sup-norm
\[\|v\|_\infty = \sup_{x\in \overline{B}_R(x_*)}|v(x)|.\]

We consider the following boundary value problem governed by a reaction-diffusion equation:
\beq
\label{e-PAppl}
\begin{cases}
u_t(t,x)=\Delta u(t,x)+\lambda g(t,x,u(t,x)), \quad (t,x)\in (0,1)\times B_R(x_*),\\
u(0,x)=\lambda \alpha(x) \beta(u(\cdot, x_*)), \quad x\in B_R(x_*),\\
u(t,x)=0, \quad (t,x)\in (0,1)\times \partial B_R(x_*),
\end{cases}
\eeq
where  $g:[0,1]\times \overline{B}_R(x_*)\times \erre \to \erre$, $\alpha:\overline{B}_R(x_*)\to \erre$,  and  $\beta: C([0,1],\erre)\to \erre$ are given non-negative functions, and $\lambda$ is a positive parameter.
\sm

Our aim is to provide the existence of a positive eigenvalue solution of the problem  \eqref{e-PAppl}. To achieve our goal, we transform the problem~\eqref{e-PAppl} into a problem of type~\eqref{e:P}. Then, we show that we can apply Theorem \ref{T-main}, under suitable assumptions on the functions $g$, $\alpha$, and $\beta$. The solution of the transformed problem assured by Theorem~\ref{T-main} will provide a corresponding solution to \eqref{e-PAppl}.
\md

\begin{theorem}\label{appthm}
Let $\rho\in (0,+\infty)$. Suppose that the following properties hold.
\begin{description}
\item{($h_g$)} The function $g:[0,1]\times \overline{B}_R(x_*)\times [0,\rho]\to \erre^+$ is continuous and such that
	\beq\label{($h_g$1)} 
	g(t,x,0)=0, \mbox{ for every } (t,x)\in [0,1]\times \partial B_R(x_*),
	\eeq
	and there exists $\mu_\rho:[0,1]\times \overline{B}_R(x_*)\to \erre^+$ for which
	\beq
	\label{($h_g$2)} 
	g(t,x,u(t,x))\ge \mu_\rho(t,x), 
	\eeq
	for every \( (t,x)\in [0,1]\times \overline{B}_R(x_*)\) and every $ u\in C([0,1]\times\overline{B}_R(x_*);\erre^+ ) $  with 
    
    $ \max_{(t,x)\in [0,1]\times\overline{B}_R(x_*) }|u(t,x)|=\rho$;
\item{($h_\alpha$)} the function $\alpha:\overline{B}_R(x_*)\to \erre^+$ is continuous and such that
	\beq\label{e-alpha0}
	\alpha(x)=0, \mbox{ for every } x\in \partial B_R(x_*);
	\eeq
\item{($h_\beta$)} the function $\beta:C([0,1];[0,\rho])\to \erre^+$
is continuous and maps bounded sets into bounded sets; moreover, assume that there exists a continuous function $\nu_\rho:\overline{B}_R(x_*)\to \erre^+$ with ${\nu_\rho}(x)=0$ if $x\in \partial B_R(x_*)$, such that 
\beq\label{e-B-eta} 
\beta(y(\cdot)(x_*))\ge \nu_\rho(x), 
\eeq
$\mbox{for every } x\in \overline{B}_R(x_*) \mbox{ and }  y\in C([0,1]; \E^+) \mbox{ with } \max_{t\in [0,1]}\|y(t)\|_{\E}=\rho.$
\end{description}
Furthermore, suppose that there exists $t_0\in\, ]0,1]$ such that
\beq
\label{e-tech} 
\max_{ x\in \overline{B}_R(x_*)}\Bigl|U(t_0)\alpha(x)\nu_\rho(x) + \int_0^{t_0} U(t_0-s) \mu_\rho(s,x)\, ds\Bigr|>0,
\eeq
where \( \{U(t)\}_{t \geq 0} \) is the $C_0$-semigroup generated by $\Delta$.

Then, 
the problem \eqref{e-PAppl} has a positive eigenvalue $\lambda_\rho$ and a positive eigenfunction $u_\rho$ such that
\( \max\left\{ |u_\rho(t,x)|:(t,x)\in [0,1]\times  \overline{B}_R(x_*) \right\}=\rho.\)
\end{theorem}

\np
{\bf Proof.}
First of all, observe that the space $C_0(B_R(x_*))$
is a Banach lattice, for the same reasons as $C(B_R(x_*))$ is  (see, e.g., \cite[Example 8.1]{AB06}). 
We put  
\beq\label{e-E}
\E:=C_0(B_R(x_*)), 
\eeq
thus the set
\[
\mathpzc{E}^+:=\{v\in \mathpzc{E}: v(x)\ge 0 \mbox{ for all } x\in B_R(x_*)\}
\]
is the positive cone of $\mathpzc{E}$. Clearly, $\mathpzc{E}^+$ is a normal cone for the  constant $c=1$ (cf. Remark \ref{r-1}).
With the same notations as in the previous section, for every $\rho\in (0,+\infty)$ we put
\[ \overline{\mathpzc{E}^+_\rho}:=\{v\in \mathpzc{E}^+: \|v\|_\mathpzc{E}\le \rho\};\qquad
\partial \mathpzc{E}^+_\rho:=\{v\in \mathpzc{E}^+: \|v\|_\mathpzc{E}= \rho\}.
\]

For any $u\in C\left( [0,1]\times \overline{B}_R(x_*) \right)$, we define the function $y_u:[0,1]\to \mathpzc{E}$ as
\beq\label{e-yu}
y_u(t)(x):=u(t,x), \quad 
x\in \overline{B}_R(x_*)).
\eeq
Note that the function $y_u$ is continuous. Indeed,
since  $u$ is continuous on the compact set $[0,1]\times \overline{B}_R(x_*)$, then it is uniformly continuous. Hence,  for every $\varepsilon>0$ there exists $\delta(\varepsilon)>0$ such that 
\beq
\label{u-uc}
|u(t,x)-u(\bar t,\bar x)|<\varepsilon,
\eeq
for every $(t,x),(\bar t,\bar x)\in [0,1]\times \overline{B}_R(x_*)$ with $|t-\bar t|<\delta(\varepsilon)$, $\|x-\bar x\|_n<\delta(\varepsilon)$ (here $\|\cdot\|_n$ is the euclidean norm in $\erre^n$). In particular, \eqref{u-uc} holds for $\bar x=x$, yielding
\[
 |u(t,x)-u(\bar t,x)|<\varepsilon, \mbox{ for every } |t-\bar t|<\delta(\varepsilon) \mbox{ and every } x\in \overline{B}_R(x_*).
\]
As a consequence, for every $t,\bar t\in [0,1]$ with $|t-\bar t|<\delta(\varepsilon)$ we have
\[
\|y_u(t)-y_u(\bar t)\|_\infty
=\sup_{x\in  \overline{B}_R(x_*))} |u(t,x)-u(\bar t,x)|\le\varepsilon,
\]
i.e., the (uniform) continuity of $y_u$ on $[0,1]$. 

\np
Moreover, we consider:
%
the linear operator
\beq\label{e-A}
Av=\Delta v 
\eeq
on the domain
\beq\label{e-D(A)}
D(A)=\{v\in C_0(B_R(x_*))
: \Delta v\in C_0(B_R(x_*)) \};
\eeq
the function $f:[0,1]\times \overline{\mathpzc{E}^+_\rho}\to \mathpzc{E}^+$ defined by
\beq\label{e-f}
f(t,v)(x)=g(t,x,v(x)), \quad x\in  \overline{B}_R(x_*);
\eeq
the operator $B:C([0,1];\overline{\mathpzc{E}^+_\rho})\to \mathpzc{E}^+$ given by
\beq\label{e-B}
B(y)(x)=\alpha (x) \beta(y(\cdot)(x_*)), \quad x\in  \overline{B}_R(x_*).
	   \eeq

By means of the positions \eqref{e-E}, \eqref{e-yu}, \eqref{e-A}-\eqref{e-B}, the system \eqref{e-PAppl} can be rewritten as 
a problem of type \eqref{e:P}, i.e.
\begin{equation}\label{e:PAA}
\begin{cases}
y_u'(t)= Ay_u(t)+\lambda f(t,y_u(t)) ,\ t \in [0,1],\\
y_u(0)=\lambda B(y_u).
\end{cases}
\end{equation}
For the sake of simplicity, from now on we write $y$ instead of $y_u$. 

Our aim is to prove that the mappings defined by \eqref{e-A}, \eqref{e-f}, and \eqref{e-B}  satisfy the hypotheses of Theorem \ref{T-main}.
Let us proceed by steps.
\sm

{\em Step 1.} We recall that the Laplace operator $\Delta$, if defined on the domain $D(A)$ given in \eqref{e-D(A)}, generates a positive and compact $C_0$-semigroup (cf. \cite[Lemma 3.1]{IM16}, see also  \cite[Theorem 7.2.5 and Lemma 7.2.1]{V03}), 
so that condition (H1) holds. 
\sm

{\em Step 2.} By means of ($h_g$), we show that the function $f$ defined by \eqref{e-f} satisfies condition (H2) on the domain $[0,1]\times \overline{\mathpzc{E}^+_\rho}$. 

First of all, the map $f$ is well defined. To prove it, let us fix $(t,v)\in [0,1]\times \overline{\mathpzc{E}^+_\rho}$. Then, the definition of $\overline{\mathpzc{E}^+_\rho}$ yields that  $v(x)\in [0, \rho]$ for every $x\in\overline{B}_R(x_*)$, so that the expression \eqref{e-f} still holds for $g$ defined on the domain $[0,1]\times \overline{B}_R(x_*)\times [0,\rho]$ (cf. ($h_g$)).
Moreover, the function
\(f(t,v)(\cdot)=g(t,\cdot,v(\cdot)) \) is continuous on $ \overline{B}_R(x_*)$, by the continuity of $v$ and $g$. Further, by \eqref{($h_g$1)} it holds that
\[
f(t,v)(x)=g(t,x,v(x))=g(t,x,0)=0,\ \mbox{ for every } x\in \partial B_R(x_*).
\] 
Hence, $f(t,v)$ belongs to $\mathpzc{E}$.
Furthermore, since  $g$ is a nonnegative function, then 
\[f(t,v)(x)\ge 0,\ \mbox{ for every } x\in \overline{B}_R(x_*),
\]
thus $f(t,v)$ belongs to $\mathpzc{E}^+$, the positive cone of $C_0(B_R(x_*))$ (cf. \eqref{e-E}).

We show now that $f$ is continuous on its domain. Indeed, let $(t_0,v_0)$ be arbitrarily fixed in $[0,1]\times \overline{\mathpzc{E}^+_\rho}$ and let $\varepsilon >0$ be fixed as well. Hypothesis ($h_g$) ensures that the map $g$ is continuous on the compact set $[0,1]\times \overline{B}_R(x_*)$, so that $g$ is uniformly continuous on its domain. Hence, there exists $\delta(\varepsilon)>0$ such that 
\beq\label{e-gc}
|g(t,x,p)-g(t_0,x,p_0)|<\varepsilon, \ \mbox{ for every $x\in \overline{B}_R(x_*)$},
\eeq
for every $t\in [0,1]$ with $|t-t_0|< \delta(\varepsilon)$,  $p,p_0\in \erre$ with $|p-p_0|<\delta(\varepsilon)$. 
Now, for any $(t,v)\in [0,1]\times \overline{\mathpzc{E}^+_\rho}$ with $|t-t_0|< \delta(\varepsilon)$ and $\|v-v_0\|_\mathpzc{E}<\delta(\varepsilon)$, we have
\[|v(x)-v_0(x)|<\delta(\varepsilon), \ \mbox{ for every } x\in \overline{B}_R(x_*),\]
so by \eqref{e-gc} we get 
\[|g(t,x,v(x))-g(t_0,x,v_0(x))|<\varepsilon, \ \mbox{ for every $x\in \overline{B}_R(x_*)$},\]
and therefore
\[
\|f(t,v)-f(t_0,v_0)\|_\mathpzc{E} =\sup_{x\in \overline{B}_R(x_*)} |g(t,x,v(x))-g(t_0,x,v_0(x))|<\varepsilon.
\]

The function $g$ is continuous on the compact set $[0,1]\times \overline{B}_R(x_*)\times [0,\rho]$, so it is bounded. Since for every $(t,v)\in [0,1]\times \overline{\mathpzc{E}^+_\rho}$ it holds that 
 \[
 \|f(t,v)\|_{\E}=\max_{x\in \overline{B}_R(x_*)} |f(t,v)(x)| =\max_{x\in \overline{B}_R(x_*)} |g(t,x,v(x))|, 
 \]
 then $f$ is bounded too.
 
Finally, let us fix any $t\in [0,1]$ and $y\in C([0,1]; \E^+)$ with $\max_{t\in [0,1]}\|y(t)\|_\E=\rho$.
It is easy to check that 
\[y(t)(x)\in [0,\rho], \mbox{ for every $x\in \overline{B}_R(x_*)$ and $t\in [0,1]$}.
\]
Hence, by \eqref{($h_g$2)} we have
\[
f(t,y(t))(x)=g(t,x,y(t)(x))\ge \mu_\rho(t,x), \  \mbox{ for every $x\in \overline{B}_R(x_*)$ and $t\in [0,1]$}.
\] 
Thus, $f$ satisfies property \eqref{e-fdelta} taking $\delta_\rho:[0,1]\to \E^+$ as  
\beq\label{e-deltamu}
\delta_\rho(t)(x):= \mu_\rho(t,x),\ \mbox{ for every } x\in \overline{B}_R(x_*).
\eeq

\sm

{\em Step 3.}  We demonstrate, by means of  ($h_\alpha$) and  ($h_\beta$), that the operator $B$ defined by \eqref{e-B} satisfies condition (H3) on $C\big([0,1];\overline{\mathpzc{E}^+_\rho}\big)$.

The operator is well posed. Indeed, fixed $y\in C\big([0,1];\overline{\E^+_\rho}\big)$, the map
\[x\mapsto B(y)(x)=\alpha (x) \beta(y(\cdot)(x_*))\]
is continuous on $\overline{B}_R(x_*)$ by the continuity of $\alpha$. Further, by \eqref{e-alpha0}, we have
\[B(y)(x)=0, \mbox{ for every } x\in \partial B_R(x_*).\]
Recalling that $\alpha,\, \beta$ are nonnegative functions, we can conclude that $B(y)\in \E^+$.

Now, we show that $B$ is continuous. Let $y,\, y_0,\in C\big([0,1];\overline{\E^+_\rho}\big)$ and consider the functions $t\mapsto y(t)(x_*)$ and $t\mapsto y_0(t)(x_*)$ defined on $[0,1]$. It is easy to see that those functions are continuous and that $y(t)(x_*),\, y_0(t)(x_*)\in [0,\rho]$, i.e. both the functions belong to $C([0,1];[0,\rho])$. Since $\beta$ is continuous on $C([0,1];[0,\rho])$ (cf. ($h_\beta$)), for every $\varepsilon>0$ there exists $\delta(\varepsilon)>0$ such that 
\[
\max_{t\in [0,1]}|y(t)(x_*)-y_0(t)(x_*)|<\delta(\varepsilon) \, \Longrightarrow \, |\beta(y(\cdot)(x_*))-\beta(y_0(\cdot)(x_*))|<\varepsilon.
\]
Hence, if 
\[\|y-y_0\|_{C\big([0,1];\overline{\E^+_\rho}\big)}:=\max_{t\in [0,1]}\max_{x\in \overline{B}_R(x_*)} |y(t)(x)-y_0(t)(x)|<\delta(\varepsilon),
\]
then
\begin{align*}
\|B(y)-B(y_0)\|_{\E}&=\max_{x\in \overline{B}_R(x_*)}|\alpha(x)| |\beta(y(\cdot)(x_*))-\beta(y_0(\cdot)(x_*))|\\
	&\le A \max_{x\in \overline{B}_R(x_*)} |\beta(y(\cdot)(x_*))-\beta(y_0(\cdot)(x_*))|
	\le A\varepsilon,
\end{align*}
where $A:=\max_{x\in \overline{B}_R(x_*)} |\alpha(x)|$.

Further, the operator $B$ is bounded, since both $\alpha$ and $\beta$ are bounded.

Finally, for every $y\in C([0,1]; \E^+)$ with $\max_{t\in [0,1]}\|y(t)\|_\E=\rho$,
by \eqref{e-B-eta} we obtain
\[
B(y)(x)\equiv\alpha(x) \beta(y(\cdot)(x_*))\ge\alpha(x)\nu_\rho(x) , \mbox{ for every } x\in \overline{B}_R(x_*).
\]
Hence, \eqref{e-Beta} is satisfied by $B$, just taking 
\beq\label{e-etanu}
\eta_\rho:=\alpha\nu_\rho.
\eeq

{\em Step 4.} By \eqref{e-tech}, we have that 
\eqref{e-H4} true for the functions $\eta_\rho$ and $\delta_\rho$ defined by \eqref{e-etanu} and \eqref{e-deltamu}, respectively.

Hence, we can apply Theorem \ref{T-main} and claim that there exists $\lambda_\rho\in (0,+\infty)$ for which the parametric nonlocal boundary value problem \eqref{e:PAA} admits a positive continuous solution  $y_\rho$ such that $\max_{t\in [0,1]}\|y_\rho(t)\|_{\E}=\rho$.
The corresponding function 
\[u_\rho(t,x):=y_\rho(t)(x), \ (t,x)\in [0,1 ]\times\overline{B}_R(x_*),
\]
satisfies the thesis.
\hfill $\Box$

In the following example we show that the bounds that occur in our theory can be either computed or estimated.
\begin{example}

\np
Let us consider the system
\beq
\label{e-Ex}
\begin{cases}
u_t(t,x)= u_{xx}(t,x)+\lambda \, tx(\pi - x)u^2(t,x), \quad (t,x)\in (0,1)\times (0,\pi),\\
u(0,x)=\lambda \, \sin x \int_0^1 e^{u(t,\pi/2)} \, dt, \quad x\in (0,\pi),\\
u(t,0)=u(t,\pi)=0, \quad t\in (0,1).
\end{cases}
\eeq
Note that the system \eqref{e-Ex} is a particular case of \eqref{e-PAppl}, for $x_*=\frac\pi2$, $R=\frac\pi2$, $\Delta=\frac{\partial^2}{\partial x^2}$, $g(t,x,p)=tx(\pi - x)p^2$, $\alpha(x)=\sin x$, and  $\beta(u(\cdot,x_*))= \int_0^1 e^{u(t,\pi/2)} \, dt$.
Note also that the conditions ($h_g$), ($h_\alpha$), and ($h_\beta$) are satisfied by these functions. In particular, we have $\mu_\rho(t,x)=0$, $(t,x)\in [0,1]\times [0,\pi]$, and $\nu_\rho(x)=1$, $x\in [0,\pi].$
As a consequence, also property~\eqref{e-tech} is satisfied by taking $t_0=1$. In fact, it holds that
\begin{multline*}
\max_{ x\in [0,\pi]}\Bigl|U(1)\sin x\int_0^1 e^{u(t,\pi/2)} \, dt+ \int_0^{t_0} U(t_0-s) \mu_\rho(s,x)\, ds\Bigr|\\
\ge \max_{ x\in [0,\pi]}\Bigl|U(1)\sin x \int_0^1 e^{u(t,\pi/2)} \, dt\Bigr|\\
     \ge \max_{ x\in [0,\pi]}|U(1)\sin x|=\max_{ x\in [0,\pi]}\Bigl|\frac{\sin x}e\Bigr|=\frac1e>0.
\end{multline*}
By a direct application of Theorem~\ref{appthm}, as in Corollary~\ref{CorEig}, we obtain the existence of uncountably many couples $(\lambda_{\rho},y_{\rho})$, of non-negative eigenvalues and non-negative eigenfuntions with localized norm, that satisfy the problem~\eqref{e-Ex}.
\end{example}

\bg

\np
{\bf Acknowledgments and funding.} 

This study was partly funded by the Unione europea - Next Generation EU, Missione 4 Componente C2 -  CUP Master: J53D2300390 0006 - Research project of MUR (Italian Ministry of University and Research) PRIN 2022  “Nonlinear differential problems with applications to real phenomena” (Grant Number: 2022ZXZTN2).

The authors are members of the ``Gruppo Nazionale per l'Analisi Matematica, la Probabilità e le loro Applicazioni'' (GNAMPA) of the Istituto Nazionale di Alta Matematica (INdAM) and of the UMI Group TAA. G. Infante is a member of the “The Research ITalian network on Approximation (RITA)”.
\md

\np
{\bf Conflicts of interest.} The authors declare no conflict of interest.
\md

\np
{\bf Contribution statement.}
Both authors contributed equally to this manuscript.
\md



\begin{thebibliography}{999}
%
\bibitem{AB06}
C. D. Aliprantis and K. C. Border, \textit{Infinite dimensional analysis. A hitchhiker's guide}, Third edition, Springer, Berlin, 2006. 
%

\bibitem{b:b-t-v}
I. Benedetti, V. Taddei and M. V{\"a}th, Evolution problems with nonlinear nonlocal boundary conditions,
\textit{J. Dynam. Differential Equations}, 
\textbf{25} (2013), 477--503.

\bibitem{b:Bou1} 
A. Boucherif, 
Nonlocal problems for parabolic inclusions, 
\textit{Discrete Contin. Dyn. Syst.}, 
\textbf{suppl.} (2009), 82--91.

\bibitem{b:Bou3} 
A. Boucherif, 
Nonlocal conditions for lower semicontinuous parabolic inclusions, 
\textit{Adv. Difference Equ.}, \textbf{2011}, Art. ID 109570, 7 pp. 

\bibitem{B91}
L. Byszewski, Theorems about the existence and uniqueness of solutions of a semilinear evolution nonlocal Cauchy problem, \textit{J. Math. Anal. Appl.}, 
\textbf{162}
(1991) 494–505.

\bibitem{CPR15}
T. Cardinali, R. Precup and P. Rubbioni, A unified existence theory for evolution equations and systems under nonlocal conditions, 
\textit{J. Math. Anal. Appl.}, 
\textbf{432} (2015), 1039--1057.
%
\bibitem{EN06}
K.-J. Engel and R. Nagel, \textit{A short course on operator semigroups}, Universitext, Springer, New York, 2006.
%

\bibitem{KOZ}
M. Kamenskii, V. Obukhovskii and P. Zecca, 
\textit{Condensing
Multivalued Maps and Semilinear Differential Inclusions in Banach
Spaces}, De Gruyter Ser. Nonlinear Anal. Appl.~7, Walter de
Gruyter, Berlin - New York, 2001.
%
\bibitem{Krasno}
M. A. Krasnosel'ski\u{i},
\textit{Positive solutions of operator equations}, Noordhoff, Groningen, 1964.

\bibitem{Kra-Lady}
M. A. Krasnosel'ski\u{\i} and L. A. Lady\v{z}enski\u{\i}, The structure of the spectrum of positive nonhomogeneous operators, \textit{Trudy Moskov. Mat. Ob\v{s}\v{c}}, \textbf{3} (1954), 321--346.

\bibitem{IM16}
G. Infante and M. Maciejewski, 
Multiple positive solutions of parabolic systems with nonlinear, nonlocal initial conditions, 
\textit{J. Lond. Math. Soc.}, \textbf{94} (2016), 859--882.
%
\bibitem{MR16}
L. Malaguti and P. Rubbioni, Nonsmooth feedback controls of nonlocal dispersal models, \textit{Nonlinearity}, 
\textbf{29} (2016), 823--850.
%
\bibitem{O13}
L. Olszowy, Existence of mild solutions for the semilinear nonlocal problem in Banach spaces, \textit{Nonlinear Anal.}, \textbf{81} (2013), 211--223.
%
\bibitem{P83}
A. Pazy,
\textit{Semigroups of linear operators and applications to partial differential equations},
Applied Mathematical Sciences, 44, Springer-Verlag, New York, 1983.
%
\bibitem{V03}
 I. I. Vrabie,
$C_0$-Semigroups and Applications,
North-Holland Mathematics Studies, 191. North-Holland Publishing Co., Amsterdam, 2003.
%
\bibitem{YuWa22}
Y.-Y. Yu and  F.-Z. Wang,  
Solvability for a nonlocal dispersal model governed by time and space integrals, 
\textit{Open Math}, \textbf{20} (2022), 1785--1799.
\end{thebibliography}
\end{document}